\documentclass{amsart}
\usepackage{psfig}    

\newtheorem{theorem}{Theorem}[section]

\newtheorem{cor}[theorem]{Corollary}
\newtheorem{conj}[theorem]{Conjecture}
\newtheorem{prop}[theorem]{Proposition}

\theoremstyle{definition}

\theoremstyle{remark}

\numberwithin{equation}{section}

\newcommand{\abs}[1]{\lvert#1\rvert}

\DeclareSymbolFont{AMSb}{U}{msb}{m}{n}
\DeclareMathSymbol{\Z}{\mathalpha}{AMSb}{"5A}

\DeclareMathSymbol{\nmid}{\mathrel}{AMSb}{"2D}

\begin{document}
\newcommand{\beqs}{\begin{equation*}}
\newcommand{\eeqs}{\end{equation*}}
\newcommand{\beq}{\begin{equation}}
\newcommand{\eeq}{\end{equation}}
\newcommand\mylabel[1]{\label{#1}}
\newcommand\eqn[1]{(\ref{eq:#1})}
\newcommand\gausspoly[2]{\begin{bmatrix} #1 \\ #2\end{bmatrix}_q}

\title[K.~Saito's Conjecture for Nonnegative Eta Products]{K.~Saito's Conjecture for 
Nonnegative Eta Products}

\author{Alexander Berkovich}
\address{Department of Mathematics, University of Florida, Gainesville,
Florida 32611-8105}
\email{alexb@math.ufl.edu}          

\author{Frank G. Garvan}
\address{Department of Mathematics, University of Florida, Gainesville,
Florida 32611-8105}
\email{frank@math.ufl.edu}          


\date{July 25, 2006}  



\begin{abstract}
We prove that the Fourier coefficients of a certain general
eta product considered by K.~Saito
are nonnegative. The proof is elementary and depends on a
multidimensional theta function identity. The $z=1$ case 
is an identity for the generating function for $p$-cores due to
Klyachko \cite{K1} and Garvan, Kim and Stanton \cite{GKS}.
\end{abstract}

\maketitle

\section{Introduction} \label{sec:intro}

Throughout this paper $q=\exp(2\pi i\tau)$ with $\Im \tau >0$ so that $\abs{q}<1$. 
As usual the
Dedekind eta function is defined as
\beq
\eta(\tau) := \exp(\pi i\tau/12) \prod_{n=1}^\infty(1- \exp(2\pi in\tau))
= q^{1/24} \prod_{n=1}^\infty(1-q^n).
\mylabel{eq:etadef} 
\eeq
An eta product is a finite product of the form
\beq
\prod_k \eta(k\tau)^{e(k)},
\mylabel{eq:etaproddef} 
\eeq
where the $e(k)$ are integers. K.~Saito \cite{S1} considered eta products
that are connected with elliptic root systems and considered the
problem of determining when all the Fourier coefficients of such eta products
are nonnegative. Subsequent work contains the following 
\begin{conj} 
\label{conj1}
{\rm (K.~Saito \cite{S3})} Let $N$ be a positive integer.
The eta product
\beq
S_N(\tau):=\frac{ \eta(N\tau)^{\phi(N)}}
{\prod_{d\mid N} \eta(d\tau)^{\mu(d)}}
\mylabel{eq:setaproddef} 
\eeq
has nonnegative Fourier coefficients. 
\end{conj}
The conjecture has been proved for $N=2,3,4,5,6,7,10$ by K.~Saito \cite{S1}, \cite{S2},
\cite{S3}, \cite{S4}, \cite{S5}, for prime powers $N=p^\alpha$ by 
T.~Ibukiyama \cite{I1}, 
and for $\gcd(N,6)>1$ by K.~Saito and S.~Yasuda \cite{SY},
who also showed that for general $N$, the coefficient of $q^n$
in $S_N(\tau)$ is nonnegative for sufficiently large $n$.
We prove the conjecture for general $N$. 

The case $N=p$ (prime) occurs in the study of $p$-cores.
A partition is a $p$-core if it has no hooks of length $p$ \cite{GKS},
\cite{JK}. 
$p$-cores are important in the study of $p$-modular
representations of the symmetric group $S_n$. Define
\beq
E(q) := \prod_{n=1}^\infty (1 - q^n),
\mylabel{eq:Edef} 
\eeq
and let $a_t(n)$ denote the number of partitions of $n$ that are
$t$-cores. It is well known that for any positive integer $t$
\beq
\sum_{n\ge0} a_t(n) q^n =  \frac{E(q^t)^t}{E(q)}.
\mylabel{eq:pcore1} 
\eeq
This result is originally due to Littlewood \cite{L}.
See \cite{GKS} for a combinatorial proof. 
Thus \eqn{pcore1} implies that Conjecture \ref{conj1} holds for
$N=p$ prime since
\beq
S_p(\tau) = \frac{\eta(p\tau)^p}{\eta(\tau)}
= q^{(p^2-1)/24} \frac{E(q^p)^p}{E(q)}.
\mylabel{eq:Sp} 
\eeq
Granville and Ono \cite{GO} have proved that $a_t(n) > 0$
for all $t\ge4$ and all $n$.
We also need the following identity due to Klyachko \cite{K1}
\beq
\sum_{\substack{\vec n\in \mathbb Z^t \\ \vec n\cdot\vec{1}_t=0}} q^{\frac{t}{2}
\vec n\cdot\vec{n}+\vec b_t\cdot\vec{n}}=
\frac{E(q^t)^t}{E(q)}, 
\mylabel{eq:kid} 
\eeq
where $\vec{1}_t=(1,1,\ldots,1)\in \mathbb Z^t$, $\vec b_t=(0,1,2,\ldots,t-1)$,
and $t$ is any positive integer.
See \cite{GKS} for a combinatorial proof. See also \cite[Prop.1.29]{G94}
and \cite[\S2]{BG06b}.
Our proof of K.~Saito's Conjecture
depends on the following extension of \eqn{kid}.
\begin{theorem}
\mylabel{thm1} Let $a\ge 2$ be an integer. Then for $z\ne0$ and $\abs{q}<1$
we have
\begin{align}
&C_a(z;q):= \mylabel{eq:Cazq}\\ 
&\sum_{\substack{\vec n=(n_0,n_1,\dots, n_{a-1})\in \mathbb Z^a \\ \vec n\cdot\vec{1}_a=0}} 
q^{\frac{a}{2}\vec n\cdot\vec{n}+\vec b_a\cdot\vec{n}}
\left( z^{an_1+1} + z^{an_2+2} + \cdots + z^{an_{a-1}+a-1} + z^{-an_{a-1}}
\right) \nonumber\\
&= E(q) E(q^a)^{a-2}
\prod_{n=1}^\infty \frac{(1-z^a q^{a(n-1)})(1-z^{-a}q^{an})}
                        {(1-z q^{n-1})(1-z^{-1}q^n)}
\nonumber
\end{align}
\end{theorem}
We note that \eqn{kid} follows from \eqn{Cazq} by letting $z\to1$.
The case $a=3$ is equivalent to \cite[(1.23)]{HGB}. See \cite[\S3.3]{BG06a}.
The case $a=2$ can be written as
\beq
\sum_{n\in\mathbb Z} q^{2n^2+n}\left(z^{2n+1}+z^{-2n}\right)
= 
\prod_{n=1}^\infty (1+z q^{(n-1)})(1+z^{-1}q^{n})(1 - q^n)
\mylabel{eq:Cazq2}, 
\eeq
which follows easily from Jacobi's triple product identity
\cite[(2.2.10)]{Andbook}.

\subsection*{Notation} 
We use the following notation for finite products
$$
(z;q)_n=(z)_n=
\begin{cases}
\prod_{j=0}^{n-1}(1-zq^j), & n>0 \\
1,                         & n=0.
\end{cases}
$$
For infinite products we use
$$
(z;q)_\infty=(z)_\infty = \lim_{n\to\infty} (z;q)_n
=\prod_{n=1}^\infty (1-z q^{(n-1)}),
$$
and
$$
[z;q]_\infty = (z;q)_\infty (z^{-1}q;q)_\infty=
\prod_{n=1}^\infty (1-z q^{(n-1)})(1-z^{-1}q^{n}),
$$
for $\abs{q}<1$ and $z\ne 0$.


\section{Proof of Theorem \ref{thm1}}
\label{sec:proofthm1}
Suppose $a\ge2$. 
The idea is to show both sides of \eqn{Cazq} satisfy the
same functional equation as $z\to zq$ and agree for enough
values of $z$.
Define
\beq
R_a(z;q) =  E(q) E(q^a)^{a-2}
            \frac{[z^a;q^a]_\infty}{[z;q]_\infty},
\mylabel{eq:Radef} 
\eeq
which is the right side of \eqn{Cazq}.
An easy calculation gives
\beq
R_a(zq;q) = z^{-(a-1)} R_a(z;q).
\mylabel{eq:Rafe} 
\eeq
We show that basically the $a$ terms in the definition of $C_a(z;q)$
are permuted cyclically as $z\to zq$. To this end we define
\beq
Q_a(\vec{n}) = \frac{a}{2}\vec n\cdot\vec{n}+\vec b_a\cdot\vec{n},
\mylabel{eq:Qadef} 
\eeq
\beq
F_j(z;q) :=
\sum_{\substack{\vec n=(n_0,n_1,\dots, n_{a-1})\in \mathbb Z^a \\ \vec n\cdot\vec{1}_a=0}} 
 z^{an_j + j} q^{Q_a(\vec{n})}
\qquad (1\le j \le a-1),
\mylabel{eq:Fjdef} 
\eeq
and
\beq
F_0(z;q) :=
\sum_{\substack{\vec n=(n_0,n_1,\dots, n_{a-1})\in \mathbb Z^a \\ 
\vec n\cdot\vec{1}_a=0}} 
 z^{-an_{a-1}} q^{Q_a(\vec{n})}.
\mylabel{eq:F0def} 
\eeq

Now suppose $2 \le j \le a-2$.
Let $\vec{e}_0=(1,0,\dots,0)$, 
$\vec{e}_1=(0,1,\dots,0)$, \dots,
$\vec{e}_{a-1}=(0,0,\dots,0,1)$ be the standard unit vectors,
$\vec{n}=(n_0,n_1,\dots, n_{a-1})\in \mathbb Z^a$,
and
$\vec{n}\,'=(n_1,n_2,\dots, n_{a-1},n_0) + \vec{e}_{j-1} - \vec{e}_{a-1}$.
An easy calculation gives
\beq
Q_a(\vec{n}\,') - Q_a(\vec{n})
= a n_j + j - \vec{n}\cdot\vec{1}_a.
\mylabel{eq:Qdiff} 
\eeq
Hence
\begin{align}
F_{j-1}(z;q) &= \sum_{\substack{\vec n\in \mathbb Z^a \\ 
\vec n\cdot\vec{1}_a=0}}
z^{an_{j-1}+(j-1)} q^{Q_a(\vec{n})}\mylabel{eq:Fjm1}\\ 
&= \sum_{\substack{\vec{n}\,'\in \mathbb Z^a \\ 
\vec n\,'\cdot\vec{1}_a=0}}
z^{a(n_{j}+1)+(j-1)} q^{Q_a(\vec{n}\,')}\nonumber\\
&= \sum_{\substack{\vec{n}\in \mathbb Z^a \\ 
\vec n\cdot\vec{1}_a=0}}
z^{an_{j}+j+(a-1)} q^{Q_a(\vec{n})+an_j+j}\nonumber\\
&=z^{(a-1)} F_j(zq;q),
\nonumber
\end{align}
and
\beq
F_j(zq;q) = z^{-(a-1)} F_{j-1}(z;q).
\mylabel{eq:Fjfe} 
\eeq
Similarly we find that
\beq
F_0(zq;q) = z^{-(a-1)} F_{a-1}(z;q),
\mylabel{eq:F0fe} 
\eeq
by using the result that
\beq
Q_a(\vec{n}\,') - Q_a(\vec{n})
= -a n_{a-1} - (a-2)\vec{n}\cdot\vec{1}_a,
\mylabel{eq:Qdiff0} 
\eeq
where
$\vec{n}\,'=(-n_{a-2},-n_{a-3},\dots, -n_1,-n_0,-n_{a-1})$.
Also we have
\beq
F_1(zq;q) = z^{-(a-1)} F_{0}(z;q),
\mylabel{eq:F1fe} 
\eeq
by using the result that
\beq
Q_a(\vec{n}\,') - Q_a(\vec{n})
= a n_1+1 - a\vec{n}\cdot\vec{1}_a,
\mylabel{eq:Qdiff1} 
\eeq
where
$\vec{n}\,'=(-n_{0},-n_{a-1},-n_{a-2},\dots, -n_2,-n_1)+\vec{e}_0 
- \vec{e}_{a-1}$.

Since
\beq
C_a(z;q) = \sum_{j=0}^{a-1} F_j(z;q),
\mylabel{eq:Casum} 
\eeq
we have
\beq
C_a(zq;q) = z^{-(a-1)} C_a(z;q).
\mylabel{eq:Cafe} 
\eeq
In view of \cite[Lemma 2]{ASD} or \cite[Lemma 1]{HGB}
it suffices to show that \eqn{Cazq} holds for $a$ distinct
values of $z$ with $\abs{q} < \abs{z} \le 1$.
It is clear that
\beq
C_a(z;q) = R_a(z;q) = 0,
\mylabel{eq:Cazqzero} 
\eeq
for $z=\exp(2\pi ik/a)$ for $1 \le k \le a-1$.
Finally, \eqn{Cazq} holds for $z=1$ since
\beq
C_a(1;q) = a 
\sum_{\substack{\vec n\in \mathbb Z^a \\ \vec n\cdot\vec{1}_a=0}} 
q^{\frac{a}{2}
\vec n\cdot\vec{n}+\vec b_a\cdot\vec{n}}
= a\frac{E(q^a)^a}{E(q)}=R_a(1;q),
\eeq
by \eqn{kid} with $t=a$.
This completes the proof of Theorem \ref{thm1}.

\section{Proof of K.~Saito's Conjecture} 
\label{sec:proofSaito}

First we show that
\beq
\prod_{d\mid M} E(q^d)^{\mu(d)} =
\prod_{\substack{ n\ge 1\\ (n,M)=1}} (1 - q^n),
\mylabel{eq:Eprop} 
\eeq
for any positive integer $M$.
Now
\beq
\prod_{d\mid M} E(q^d)^{\mu(d)}
=
\prod_{d\mid M} \prod_{m=1}^\infty (1- q^{dm})^{\mu(d)}
=
\prod_{n=1}^\infty (1- q^{n})^{\varepsilon(n)},
\mylabel{eq:Eep} 
\eeq
where
\beq
\varepsilon(n) =
\sum_{d\mid M\,\&\,d\mid n} \mu(d)
=
\sum_{d\mid (M,n)} \mu(d)
=
\begin{cases}
1 & \mbox{if $(M,n)=1$}\\
0 & \mbox{otherwise},
\end{cases}
\mylabel{eq:epsimp} 
\eeq
by a well known property of the M\"obius function, and we have
\eqn{Eprop}.

For any positive integer $N$ we define
\beq
S_N(q) := \frac
{E(q^N)^{\phi(N)}}
{\prod_{d\mid N} E(q^d)^{\mu(d)}}.
\mylabel{eq:SNdef} 
\eeq
We wish to show that all coefficients in the $q$-expansion of $S_N(q)$
are nonnegative. We consider three cases.

\noindent
{\it Case 1.} $N=p^\alpha$ where $p$ is prime. This case was proved
by Ibukiyama \cite{I1}. Alternatively, the case $\alpha=1$ follows from 
\eqn{pcore1} and then use an easy induction on $\alpha$.

\noindent
{\it Case 2.} $N=p\,M$, where $p$ is prime, $M$ is odd and $p\nmid M$.
We have
\beq
\prod_{d\mid N} E(q^d)^{\mu(d)} =
\prod_{d\mid M} E(q^d)^{\mu(d)}  E(q^{pd})^{\mu(pd)} =
\prod_{d\mid M} \left(
\frac{E(q^d)}{E(q^{pd})}\right)^{\mu(d)}.
\mylabel{eq:Eprod} 
\eeq
By \eqn{Eprop} we have
\beq
\prod_{d\mid M} E(q^d)^{\mu(d)} =
\prod_{\substack{ n\ge 1\\ (n,M)=1}} (1 - q^n) =
\prod_{n\ge0}
\prod_{\substack{(r,M)=1\\1\le r \le M-1}} (1 - q^{Mn+r}) =
\prod_{\substack{(r,M)=1\\1\le r \le \frac{M-1}{2}}} [q^r;q^M]_\infty.
\mylabel{eq:Eprod2} 
\eeq
Now for $a$ a positive integer, $\abs{q}<1$ and $z\ne0$ we let
\beq
D_a(z;q) := \frac{E(q^a)^a}{E(q)} C_a(z;q)
         = E(q^a)^{2a-2} \frac{[z^a;q^a]_\infty}{[z;q]_\infty}
\mylabel{eq:Dazq}, 
\eeq
so that
\begin{align}
\prod_{\substack{(r,M)=1\\1\le r \le \frac{M-1}{2}}} D_p(q^r;q^M) 
&= \left(E(q^{pM})^{2p-2}\right)^{\phi(M)/2}
\prod_{\substack{(r,M)=1\\1\le r \le \frac{M-1}{2}}} 
\frac{[q^{pr};q^{pM}]_\infty}{[q^{r};q^{M}]_\infty} \nonumber\\
&=E(q^N)^{\phi(N)}
\prod_{d\mid M} \frac{ E(q^{pd})^{\mu(d)}}{ E(q^{d})^{\mu(d)}}
\qquad\mbox{(by \eqn{Eprod2})}
\mylabel{eq:Dprod}\\ 
&= \frac
{E(q^N)^{\phi(N)}}
{\prod_{d\mid N} E(q^d)^{\mu(d)}} 
\qquad\mbox{(by \eqn{Eprod})}
\nonumber\\
&= S_N(q).
\nonumber
\end{align}
K.~Saito's Conjecture holds in this case since each $C_p(q^r;q^M)$
has nonnegative coefficients by Theorem \ref{thm1}, and
$E(q^{pM})^p/E(q^M)$ has nonnegative coefficients by \eqn{pcore1}
so that each $D_p(q^r;q^M)$ has nonnegative coefficients.

\noindent
{\it Case 3.} $N=p^\alpha\,M$, where $p$ is prime, $M$ is odd, $p\nmid M$,
and $\alpha\ge2$. We let $N'=p\,M$. It is clear that
\beq
\prod_{d\mid N} E(q^d)^{\mu(d)}= \prod_{d\mid N'} E(q^d)^{\mu(d)}.
\mylabel{eq:Eprop2} 
\eeq
Hence                              
\beq
S_N(q) = \frac{E(q^N)^{\phi(N)}}
{E(q^{N'})^{\phi(N')}} S_{N'}(q) 
= \left(\frac{E(q^{p^{\alpha-1}N'})^{p^{\alpha-1}}}{E(q^{N'})}
\right)^{ (p-1)\phi(M)} S_{N'}(q).
\mylabel{eq:Sprop} 
\eeq
Here $S_N(q)$ is the product of two terms. The second term
$S_{N'}(q)$ has 
nonnegative coefficients from Case 2. The first term has
nonnegative coefficients
using \eqn{pcore1} with $q$ replaced
with $q^{N'}$ and $t=p^{\alpha-1}$. Thus K.~Saito's Conjecture holds
in this case.

\section{Other Products with Nonnegative Coefficients} 
\label{sec:otherprods}

In this section we state a number results for coefficients
of other infinite products. More detail will appear in a later version
of this paper. For a formal power series
$$
F(q) := \sum_{n=0}^\infty a_n q^n \in \Z[[q]]
$$
we write
$$
F(q) \succeq 0,
$$
if
$a_n \ge 0$ for all $n\ge0$.
For a formal power series $F(z_1,z_2,\dots,z_n;q)$ in more than one variable
we interpret $F(z_1,z_2,\dots,z_n;q)\succeq0$ in the natural way.
The following result follows from the $q$-binomial theorem
\cite[Thm2.1]{Andbook}.
\begin{prop}
\label{atineq}If $\abs{q}$, $\abs{t} < 1$ then
\beq
\frac{ (at;q)_\infty}{(a;q)_\infty (t;q)_\infty}
=\sum_{n=0}^\infty 
\frac{t^n}{(aq^n;q)_\infty (q)_n} \succeq 0.
\mylabel{eq:atq} 
\eeq
\end{prop}

\begin{cor}
\label{atcor}

\vskip 0pt
\begin{enumerate}
\item[(i)]
If $a$, $b$, $M\ge1$ then
\beq
\prod_{n=0}^\infty
\frac{ (1 - q^{Mn+a+b}) }{ (1-q^{Mn+a}) (1-q^{Mn+b}) } \succeq 0.
\mylabel{eq:coratq1} 
\eeq
\item[(ii)]
If $m$, $n>1$ then
\beq
\frac{ E(q^m) E(q^n) }{E(q)} \succeq0.
\mylabel{coratq2}
\eeq
\item[(iii)]
If $m$, $n>1$ and not both $2$ then
\beq
\frac{ E(q^m) E(q^n) E(q^{mn})}{E(q)} \succeq0.
\mylabel{coratq3}
\eeq
\end{enumerate}
\end{cor}

Proposition \ref{atineq} has a finite analogue. 
For $0 \le m \le n$ the Gaussian polynomial 
\cite[p.33]{Andbook} is defined by
\beq
\gausspoly{n+m}{m} =
\frac{ (q)_{m+n}}{(q)_n (q)_m}=
\frac{ (1-q^{n+1}) \cdots (1-q^{n+m})}{(q)_m}
\mylabel{eq:gpdef} 
\eeq
Since it is 
is the generating function for partitions with at most $m$ parts each $\le n$
it is a polynomial (in $q$) with positive integer coefficients.
We have
\begin{prop}
\label{atfinite}
If $L\ge0$ then
\beq
\frac{ (z_1 z_2;q)_L }{ (z_1;q)_L (z_2;q)_L }
=
\sum_{j=0}^L
\gausspoly{L}{j} \frac{ z_1^j }
{(z_1q^{L-j};q)_j (z_2q^{j};q)_{L-j}} \succeq 0.
\mylabel{eq:atqfin} 
\eeq
\end{prop}

The case $t=a^{-1}$ of Proposition \ref{atineq} is related to the
crank of partitions \cite{AG}.
Let $M(m,n)$ denote the number of partitions of $n$ with crank $m$.
Then
\begin{align}
(1-z) \frac{E(q)}{[z;q]_\infty}
&= \prod_{n=1}^\infty \frac{(1 - q^n)}{(1-zq^n)(1-z^{-1}q^n)} 
\mylabel{eq:crankgen} 
\\
&= 1 + (z + z^{-1} - 1)q + 
\sum_{n\ge2}\left(\sum_{m=-n}^n M(m,n) z^m\right) q^n.
\nonumber                       
\end{align}
We note the coefficients on the right side of \eqn{crankgen} are
nonnegative except for the coefficient of $z^0q^1$.
By observing that
$$
(1+z)(z + z^{-1} - 1) = z^2 + z^{-1}
$$
we have
\begin{prop}
\label{aci}
If $\abs{q}<1$ and $z\ne0$ then
\beq
(1-z^2) \frac{E(q)}{[z;q]_\infty}
= (1+z) \prod_{n=1}^\infty \frac{(1 - q^n)}{(1-zq^n)(1-z^{-1}q^n)}
\succeq 0.
\mylabel{eq:aci} 
\eeq
\end{prop}
\begin{prop}
\label{res}
If $\abs{q}<1$ and $z\ne0$ then
\beq
\frac{E(q^2)[z^4;q^2]_\infty}{[z^2;q^2]_\infty [qz^3;q^2]_\infty}
\succeq 0.
\mylabel{eq:res1} 
\eeq
\beq
\frac{E(q^3) (z^2;q^3)_\infty}{ (q^3 z^{-1};q^3)_\infty (z;q)_\infty)}
\succeq 0.
\mylabel{eq:res2} 
\eeq
\end{prop}


We make the following
\begin{conj}
\label{conj2}
Suppose $\abs{q}<1$ and $z\ne0$.
\begin{enumerate}
\item[(i)]
If $p\ge1$ then
\beq
\frac{ E(q) } {(z;q)_\infty (qz^{-p};q)_\infty} \succeq 0.
\mylabel{eq:conj2a} 
\eeq
\item[(ii)]
If $a$, $b$, $m$, $n\ge1$ then
\beq
\frac{E(q^{ma + nb})}
{(q^a;q^{ma+nb})_\infty (q^b;q^{ma+nb})_\infty} \succeq 0.
\mylabel{eq:conj2b} 
\eeq
\item[(iii)]
For $a\ge1$                        
\beq
\frac{[z^a;q]_\infty E(q)}
{[z;q]_\infty [z^{a+1};q]_\infty} \succeq 0.
\mylabel{eq:conj2c} 
\eeq
\end{enumerate}
\end{conj}

The case $p=1$ of \eqn{conj2a}, the case $m=n=1$ of \eqn{conj2b} 
and the case $a=1$ of \eqn{conj2a} are all special cases
of Proposition \ref{atineq}.
The case $a=2$ of \eqn{conj2c} follows from Proposition \ref{atineq}
together the following identity due to Ekin \cite[(42)]{E1}
\beq
\frac{[z^2;q]_\infty E(q)}
{[z;q]_\infty [z^{3};q]_\infty}  =
\frac{E(q^3)}{[z^3;q^3]_\infty [z^{-3}q;q^3]_\infty}
+z\frac{E(q^3)}{[z^3;q^3]_\infty [z^{-3}q^2;q^3]_\infty}
\mylabel{eq:ekin} 
\eeq
This identity was used by Ekin to prove an number of inequalities
for the crank of partitions mod $7$ and $11$.

\noindent
\textbf{Acknowledgement}

\noindent
We would like to thank Professor Kyoji Saito for his interest,
comments and for sending copies of his unpublished work \cite{S5},
\cite{SY}.


\bibliographystyle{amsplain}

\begin{thebibliography}{10}
\bibitem{Andbook}
G.~E.~Andrews,
 \emph{ The Theory of Partitions},
Encyclopedia of Mathematics and Its Applications, Vol.~2
(G.~-~C.~Rota, ed.), Addison-Wesley, Reading, Mass., 1976. (Reissued:
Cambridge Univ. Press, London and New York, 1985)
\bibitem{AG}
G. E. Andrews and F. G. Garvan,
\textit{Dyson's crank of a partition},
{Bull. Amer. Math. Soc. (N.S.)}
\textbf{18} (1988), 167--171.
\bibitem{ASD}
A.~O.~L.~Atkin and P.~Swinnerton-Dyer,
\textit{Some properties of partitions}, 
Proc. London Math. Soc.
\textbf{4} (1954), 84--106.
\bibitem{BG06a}
A. Berkovich and F.G. Garvan,
\textit{On the {A}ndrews-{S}tanley refinement of {R}amanujan's
              partition congruence modulo 5 and generalizations},
{Trans. Amer. Math. Soc.},
\textbf{358} (2006), {703--726}.                        
\bibitem{BG06b}
A. Berkovich and F.G. Garvan,
\textit{The BG-rank of a partition and its applications},
arXiv:math.CO/0602362.
\bibitem{E1}
A.B. Ekin,
\textit{Inequalities for the crank},
{J. Combin. Theory Ser. A},
\textbf{83} (1998), {283--289}.
\bibitem{GKS}
F. Garvan, D. Kim and D. Stanton,
\textit{Cranks and $t$-cores},
{Invent. Math.}
\textbf{101} (1990), 1--17.
\bibitem{G94}
F. Garvan,
\textit{Cubic modular identities of {R}amanujan, hypergeometric
              functions and analogues of the arithmetic-geometric mean
              iteration}, in
{``The Rademacher legacy to mathematics'' (University Park, PA,
              1992)},
{Contemp. Math.}, Vol.{166}, {pp. 245--264}, {Amer. Math. Soc.},
{Providence, RI}, {1994}.
\bibitem{GO}
A. Granville and K. Ono,
\textit{Defect zero {$p$}-blocks for finite simple groups},
{Trans. Amer. Math. Soc.},
\textbf{348} (1996),
{331--347}.
\bibitem{HGB}
M. Hirschhorn, F. Garvan and J.Borwein,
\textit{Cubic analogues of the Jacobian theta function $\theta(z,q)$},
{Canad. J. Math.}
\textbf{45} (1993), 673--694.  
\bibitem{I1}
T. Ibukiyama,
\textit{Positivity of eta products---a certain case of {K}. {S}aito's
              conjecture},
{Publ. Res. Inst. Math. Sci.},
\textbf{41} (2005), {683--693.
\bibitem{JK}
G. James and A. Kerber,
\textit{The Representation Theory of the Symmetric Group},
Addison-Wesley, Reading, MA, 1981.
\bibitem{K1}
A. A. Klyachko,
\textit{Modular forms and representations of symmetric groups},
Jour. Soviet Math.
\textbf{26} (1984), 1879--1887.
\bibitem{L}
D. E. Littlewood, 
\textit{Modular representations of symmetric groups},
{Proc. Roy. Soc. London. Ser. A.}
\textbf{209} (1951), 333--353.
\bibitem{S1}
K.~Saito,
\textit{Extended affine root systems. {V}. {E}lliptic eta-products and
              their {D}irichlet series}, in
``Proceedings on Moonshine and related topics'' (Montr\'eal, QC,
              1999)}, {CRM Proc. Lecture Notes}, Vol. {30},
{pp. 185--222}, {Amer. Math. Soc.}, {Providence, RI}, {2001}.
\bibitem{S2}
K.~Saito,
\textit{Duality for regular systems of weights}, 
{Asian J. Math.},
\textbf{2} (1998), 983--1047.
\bibitem{S3}
K.~Saito,
\textit{Nonnegativity of Fourier coefficients of eta-products} (Japanese),
in ``Proceedings of the Second Spring Conference on Automorphic Forms
and Related Subjects,'' Careac Hamamatsu, Feb. 2003.
\bibitem{S4}
K.~Saito,
\textit{Eta-product $\eta(7\tau)^7/\eta(\tau)$}, 
arXiv:math.NT/0602367
\bibitem{S5}
K.~Saito,
\textit{Eta-product $\eta_{\Phi_h}(\tau)$}, 
unpublished note.
\bibitem{SY}
K.~Saito and S.~Yasuda,
\textit{Non-negativity of the Fourier coefficients of certain eta products}, 
Notes, June, 2006.
%
\end{thebibliography}

\end{document}